\documentclass{amsart}

\usepackage{amsmath,amssymb,amsfonts,enumerate,amsthm, amscd}
\usepackage[french, ngerman, italian, english]{babel}
\usepackage[colorlinks=true]{hyperref}
\usepackage{paralist}
\usepackage{verbatim}
\usepackage{color}
\usepackage{xcolor}
\usepackage{tikz}
\usetikzlibrary{arrows,decorations.pathmorphing,backgrounds,positioning,fit}
\RequirePackage{pgfrcs}

\newtheorem{thm}{\bf Theorem}[section]

\newtheorem{prop}[thm]{\bf Proposition}
\newtheorem{lem}[thm]{\bf Lemma}
\newtheorem{cor}[thm]{\bf Corollary}

\theoremstyle{definition}

\newtheorem{defn}[thm]{Definition}
\theoremstyle{remark}
\newtheorem{rem}[thm]{Remark}
\newtheorem{exam}[thm]{Example}

\theoremstyle{example}

\def \Supp{{\mathrm{Supp}}}

\def \k{{[k]}}

\def \a{\mathbf a}

\def \1{\mathbf 1}
\def \h{\mathrm{ht}}

\def \depth{\mathrm{depth}}
\def \sdepth{\mathrm{sdepth}}

\def \lk{\mathrm{lk}}
\def \dl{\mathrm{dl}}
\def \Del{\Delta}

\def \NN{\mathbb N}
\def \ZZ{\mathbb Z}

\def \F{\mathcal F}

\def \D{\mathcal D}

\begin{document}

\title[$k$-shellable simplicial complexes and graphs]{$k$-shellable simplicial complexes and graphs}
\author{Rahim Rahmati-Asghar}

\address{Department of Mathematics, Faculty of Basic Sciences\\
University of Maragheh\\
Maragheh, P.O. Box: 55181-83111\\
Iran}
\email{rahmatiasghar.r@gmail.com}

\begin{abstract}
In this paper we show that a $k$-shellable simplicial complex is the expansion of a shellable complex. We prove that the face ring of a pure $k$-shellable simplicial complex satisfies the Stanley conjecture. In this way, by applying expansion functor to the face ring of a given pure shellable complex, we construct a large class of rings satisfying the Stanley conjecture.

Also, by presenting some characterizations of $k$-shellable graphs, we extend some results due to Castrill\'{o}n-Cruz, Cruz-Estrada and Van Tuyl-Villareal.
\end{abstract}

\maketitle

\section*{Introduction}

Let $\Del$ be a simplicial complex on the vertex set $X:=\{x_1,\ldots,x_n\}$. Denote by $\langle F_1,\ldots,F_r\rangle$ the simplicial complex $\Del$ with facets $F_1,\ldots,F_r$. $\Del$ is called \emph{shellable} if its facets can be given a linear order $F_1,\ldots,F_r$, called a \emph{shelling order}, such that for all $2\leq j$, the subcomplex $\langle F_j\rangle\cap\langle F_1,\ldots,F_{j-1}\rangle$ is pure of dimension $\dim(F_j)-1$ (see \cite{BrMa} for probably the earliest definition of this term or as well \cite{BjWa} for a more recent exposition). Studying combinatorial properties of shellable simplicial complexes and algebraic constructions of their face rings and also the edge ideals associated to shellable graphs is a current trend in combinatorics and commutative algebra. See for example \cite{BjWa, BrMa, CaCr, Dr, HePo, VaVi}.

In this paper, we recall from \cite{Ra} the concept of $k$-shellability, and extend some results obtained previously by researchers. Actually, $k$ is a positive integer and for $k=1$, $1$-shellability coincides with shellability.

Richard Stanley \cite{St}, in his famous article ``Linear Diophantine equations and local cohomology'', made a
striking conjecture predicting an upper bound for the depth of a multigraded module. This conjecture is nowadays
called the Stanley conjecture and the conjectured upper bound is called the Stanley depth of a module. The Stanley conjecture has become quite popular, with numerous publications dealing with different aspects of the Stanley depth. Although a counterexample has apparently recently been found to the Stanley conjecture (see \cite{DuGoKlMa}), this makes it perhaps even more interesting to explore the relationship between depth and Stanley depth.

Let $S=K[x_1,\ldots,x_n]$ be the polynomial ring over a field $K$. Dress proved in \cite{Dr} that the simplicial complex $\Del$ is shellable if and only if its face ring is clean. It is also known that cleanness implies pretty cleanness. Furthermore, Herzog and Popescu \cite[Theorem 6.5]{HePo} proved that, if $I\subset S$ is a monomial ideal, and $S/I$ is a multigraded pretty clean ring, then the Stanley conjecture holds for $S/I$. It follows that, for a shellable simplicial complex $\Del$, the face ring $K[\Del]=S/I_\Del$ satisfies the Stanley conjecture where $I_\Del$ denotes the Stanley-Reisner ideal of $\Del$. We extend this result, in pure case, by showing that the face ring of a $k$-shellable simplicial complex satisfies the Stanley conjecture (see Theorem \ref{k-sh stanley}). We obtain this result by extending Proposition 8.2 of \cite{HePo} and presenting a filtration for the face ring of a $k$-shellable simplicial complex in Theorem \ref{k-sh and k-cl}.

A simple graph $G$ is called shellable if its independence complex $\Del_G$ is a shellable simplicial complex. Shellable graphs were studied by several researchers in recent years. For example, Van Tuyl and Villarreal in \cite{VaVi} classified all of shellable bipartite graphs. Also, Castrill\'{o}n and Cruz characterized the shellable graphs and clutters by using the properties of simplicial vertices, shedding vertices and shedding faces (\cite{CaCr}).

Here, we present some characterizations of $k$-shellable graphs and extend some results of \cite{CaCr}, \cite{CrEs} and \cite{VaVi} (see Theorems \ref{G to lk}, \ref{sh iff lk} and \ref{k-sh equi}). Our idea is to define a new notion, called a $k$-simplicial set, which is a generalization of the notion of simplicial vertex defined in \cite{Di} or \cite{LeBo}.

\section{Preliminaries}

For basic definitions and general facts on simplicial complexes, we refer to Stanley's book \cite{St1}.

A simplicial complex $\Del$ is \emph{pure} if all of its facets (maximal faces) are of
the same dimension. The \emph{link} and \emph{deletion} of a face $F$ in $\Del$ are defined respectively
$$\lk_\Del(F)=\{G\in\Del:G\cap F=\emptyset\ \mbox{and}\ G\cup F\in\Del\}$$
and
$$\dl_\Del(F)=\{G\in\Del:F\nsubseteq G\}.$$

Let $G$ be a simple (no loops or multiple edges) undirected graph on the vertex set $V(G)=X$ and the edge set $E(G)$. The \emph{independence complex} of $G$ is denoted by $\Del_G$ and $F$ is a face of $\Del_G$ if and only if there is no edge of $G$ joining any two vertices of $F$. The \emph{edge ideal} of $G$ is defined a quadratic squarefree monomial ideal $I(G)=(x_ix_j:x_ix_j\in E(G))$. It is known that $I(G)=I_{\Del_G}$. We say $G$ is a shellable graph if $\Del_G$ is a shellable simplicial complex.

In the following we recall the concept of expansion functor in a combinatorial and an algebraic setting from \cite{Ra1} and \cite{BaHe}, respectively.

Let $\alpha=(k_1,\ldots,k_n)$ be an $n$-tuple with positive integer entries in $\NN^n$. For $F=\{x_{i_1},\ldots,x_{i_r}\}\subseteq X$ define
$$F^\alpha=\{x_{i_11},\ldots,x_{i_1k_{i_1}},\ldots,x_{i_r1},\ldots,x_{i_rk_{i_r}}\}$$
as a subset of $X^\alpha:=\{x_{11},\ldots,x_{1k_1},\ldots,x_{n1},\ldots,x_{nk_n}\}$. $F^\alpha$ is called \emph{the
expansion of $F$ with respect to $\alpha$.}

For a simplicial complex $\Del=\langle F_1,\ldots,F_r\rangle$ on $X$, we define \emph{the expansion of $\Del$ with
respect to $\alpha$} as the simplicial complex $\Del^\alpha=\langle F^\alpha_1,\ldots,F^\alpha_r\rangle$ (see \cite{Ra1}).

In \cite{BaHe} Bayati and Herzog defined the expansion functor in the category of finitely generated multigraded
$S$-modules and studied some homological behaviors of this functor. We recall the expansion functor defined by them
only in the category of monomial ideals and refer the reader to \cite{BaHe} for more general case in the category of
finitely generated multigraded $S$-modules.

Set $S^\alpha$ a polynomial ring over $K$ in the variables
$$x_{11},\ldots,x_{1k_1},\ldots,x_{n1},\ldots,x_{nk_n}.$$
Whenever $I\subset S$ is a monomial ideal minimally generated by $u_1,\ldots,u_r$, the expansion of $I$ with respect to
$\alpha$ is defined
$$I^\alpha=\sum^r_{i=1}P^{\nu_1(u_i)}_1\ldots P^{\nu_n(u_i)}_n\subset S^\alpha$$
where $P_j=(x_{j1},\ldots,x_{jk})$ is a prime ideal of $S^\alpha$ and $\nu_j(u_i)$ is the exponent of $x_j$ in $u_i$.

\begin{exam}
Let $I\subset K[x_1,\ldots,x_3]$ be a monomial ideal minimally generated by $G(I)=\{x^2_1x_2,x_1x_3,x_2x^2_3\}$ and let $\alpha=(2,2,1)\in\NN^3$. Then
$$\begin{array}{rl}
    I^\alpha=& (x_{11},x_{12})^2(x_{21},x_{22})+(x_{11},x_{12})(x_{31})+(x_{21},x_{22})(x_{31})^2 \\
            =& (x^2_{11}x_{21},x_{11}x_{12}x_{21},x^2_{12}x_{21},x^2_{11}x_{22},x_{11}x_{12}x_{22},x^2_{12}x_{22},x_{11}x_{31},x_{12}x_{31},x_{21}x^2_{31},x_{22}x^2_{31})
  \end{array}
$$
\end{exam}

It was shown in \cite{BaHe} that the expansion functor is exact and so $(S/I)^\alpha=S^\alpha/I^\alpha$. The following lemma implies that two above concepts of expansion functor are related.

\begin{lem}\label{epansion s-R}
(\cite[Lemma 2.1]{Ra1}) For a simplicial complex $\Del$ and $\alpha\in\NN^n$ we have $(I_\Del)^\alpha=I_{\Del^\alpha}$. In particular, $K[\Del]^\alpha=K[\Del^\alpha]$.
\end{lem}

In this paper we just study the functors $\alpha=(k_1,\ldots,k_n)\in\NN^n$ with $k_i=k_j$ for all $i,j$. For
convenience, we set $\alpha=\k$ when every component of $\alpha$ is equal to $k\in\NN$. We call $I^\k$ (resp. $\Del^\k$) the expansion of $I$ (resp. $\Del$) with respect to $k$.

\section{Some combinatorial and algebraic properties of $k$-shellable complexes}

The notion of $k$-shellable simplicial complexes was first introduced by Emtander, Mohammadi and Moradi \cite{EMM} to provide a natural generalization of shellability. It was shown in \cite[Theorem 6.8]{EMM} that a simplicial complex $\Del$ is $k$-shellable if and only if the Stanley-Reisner ideal of its Alexander dual has $k$-quotionts, i.e. there exists an ordering $u_1,\ldots,u_r$ of the minimal generators of $I_{\Del^\vee}$ such that if we for $s=1,\ldots,t$, put $I_s=(u_1,\ldots,u_s)$, then for every $s$ there are monomials
$v_{s_i}$ , $i=1,\ldots,r_s$, $\deg(v_{s_i})=k$ for all $i$, such that $I_s:u_s=(v_{s_1},\ldots,v_{s_{r_s}}).$

In \cite{Ra}, we gave another definition of $k$-shellability and having $k$-quotients by adding a condition to Emtander, Mohammadi and Moradi's. In our definition the colon ideals $I_s:u_s$ were generated by regular sequences for all $s$ and in this way, all of structural properties of monomial ideals with linear quotients were generalized. The reader is referred to \cite{HeTa} for the definition of monomial ideals with linear quotients.

\begin{defn}\label{Def k-shell}
(\cite{Ra}) Let $\Del$ be a $d$-dimensional simplicial complex on $X$ and let $k$ be an integer with $1\leq k\leq d+1$. $\Del$ is called \emph{$k$-shellable} if its facets can be ordered $F_1,\ldots,F_r$, called \emph{$k$-shelling order}, such that for all
$j=2,\ldots,r$, the subcomplex $\Del_j=\langle F_j\rangle\cap\langle F_1,\ldots,F_{j-1}\rangle$ satisfies the following
properties:
\begin{enumerate}[\upshape (i)]
  \item It is generated by a nonempty set of maximal proper faces of $\langle F_j\rangle$ of dimension $|F_j|-k-1$;
  \item If $\Del_j$ has more than one facet then for every two disjoint facets $\sigma,\tau\in\langle F_j\rangle\cap\langle F_1,\ldots,F_{j-1}\rangle$ we have
      $F_j\subseteq\sigma\cup\tau$.
\end{enumerate}
\end{defn}

\begin{rem}
It follows from the definition that two concepts 1-shellability and shellability coincide.
\end{rem}

\begin{rem}
Note that the notions of $1$-shellability in our sense and Emtander, Mohammadi and Moradi's coincide. Although, for $k>1$, a simplicial complex may be $k$-shellable in their concept and not in ours. For example, consider the complex $\Del=\langle abc,aef,cdf\rangle$ on $\{a,b,\ldots,f\}$. It is easy to check that $\Del$ is $2$-shellable in the sense of \cite{EMM} but not in ours.
\end{rem}

In the following proposition we describe some the combinatorial properties of $k$-shellable complexes.

\begin{prop}\label{k-shell}
Let $\Del$ be a $d$-dimensional (not necessarily pure) simplicial complex on $X$ and let $k$ be an integer with
$1\leq k\leq d+1$. Suppose that the facets of $\Del$ can be ordered $F_1,\ldots,F_r$. Then the following conditions are
equivalent:
\begin{enumerate}[\upshape (a)]
  \item $F_1,\ldots,F_r$ is a $k$-shelling of $\Del$;
  \item for every $1\leq j\leq r$ there exist the subsets $E_1,\ldots,E_t$ of $X$ such that the $E_i$ are mutually disjoint and $|E_i|=k$ for all $i$ and the set of the minimal elements of $\langle F_1,\ldots,F_j\rangle\backslash\langle
      F_1,\ldots,F_{j-1}\rangle$ is $\{\{a_1,\ldots,a_t\}:a_i\in E_i\ \mbox{for all}\ i\}$;

  \item for all $i,j$, $1\leq i<j\leq r$, there exist $x_1,\ldots,x_k\in F_j\backslash F_i$ and some $l\in\{1,\ldots,j-1\}$
      with $F_j\backslash F_l=\{x_1,\ldots,x_k\}$.
\end{enumerate}
\end{prop}
\begin{proof}
(a)$\Rightarrow$(b): Let
$\langle F_j\rangle\cap\langle F_1,\ldots,F_{j-1}\rangle=\langle F_j\backslash\sigma_1,\ldots,F_j\backslash\sigma_t\rangle$
where $|\sigma_i|=k$ for all $i$. Since for all $i\neq i'$, $F_j\subseteq (F_j\backslash\sigma_i)\cup(F_j\backslash\sigma_{i'})$, we have $\sigma_i\cap\sigma_{i'}=\emptyset$. Hence the minimal elements of $\langle F_1,\ldots,F_j\rangle\backslash\langle F_1,\ldots,F_{j-1}\rangle$ are in the form $\{a_1,\ldots,a_t\}$ where $a_i\in\sigma_i$ for all $i$.

(b)$\Rightarrow$(c): For all $i$, suppose that $E_i=\{x_{i1},\ldots,x_{ik}\}$. Let $1\leq i<j\leq r$ and let $\{x_{1i_1},\ldots,x_{ti_t}\}$ be a minimal element of $\langle F_1,\ldots,F_j\rangle\backslash\langle F_1,\ldots,F_{j-1}\rangle$. Because $\{x_{1i_1},\ldots,x_{ti_t}\}\nsubseteq F_i$ we may assume that $x_{1i_1}\in F_j\backslash F_i$. We claim that $x_{11},\ldots,x_{1k}\in F_j\backslash F_i$. Suppose, on the contrary, that for some $s$, $x_{1s}\not\in F_j\backslash F_i$ then $x_{1s}\in F_i$ and so $x_{1s}\not\in \langle F_1,\ldots,F_j\rangle\backslash\langle F_1,\ldots,F_{j-1}\rangle$. It follows that $\{x_{1s},x_{2i_2},\ldots,x_{ti_t}\}$ is not a minimal element of $\langle F_1,\ldots,F_j\rangle\backslash\langle F_1,\ldots,F_{j-1}\rangle$, a contradiction. Therefore $E_1\subseteq F_j\backslash F_i$.

Now suppose that for all $l<j$, if $E_1$ is contained in $F_j\backslash F_l$ then $E_1\subsetneqq F_j\backslash F_l$. Then there exists $y\in F_j\backslash E_1$ such that $\{y,x_{2i_2},\ldots,x_{ti_t}\}$ is a minimal element of $\langle F_1,\ldots,F_j\rangle\backslash\langle F_1,\ldots,F_{j-1}\rangle$ different from the elements of $\{\{a_1,\ldots,a_t\}:a_i\in E_i\ \mbox{for all}\ i\}$, a contradiction. Therefore there exists $l<j$ with $F_j\backslash F_l=E_1$.

(c)$\Rightarrow$(a): Let $F\in\langle F_j\rangle\cap\langle F_1,\ldots,F_{j-1}\rangle$. Then $F\subseteq F_i$ for some $i<j$. By the condition (c), there exist $x_1,\ldots,x_k\in F_j\backslash F_i$ and some $l\in\{1,\ldots,j-1\}$ with $F_j\backslash F_l=\{x_1,\ldots,x_k\}$. But $F_j\backslash\{x_1,\ldots,x_k\}$ is a proper face of $\langle F_j\rangle\cap\langle F_1,\ldots,F_{j-1}\rangle$, because $F_j\backslash\{x_1,\ldots,x_k\}=F_j\cap F_l$. Moreover, $F_j\backslash\{x_1,\ldots,x_k\}$ is a maximal face. Finally, since $F$ is contained in $F_j\backslash\{x_1,\ldots,x_k\}$, the assertion is completed.
\end{proof}

\begin{exam}
The Figure \ref{shell exams} indicates the pure shellable and pure 2-shellable simplicial complexes of dimensions 1, 2 and 3 with
3 facets.

\begin{figure}
$$\begin{array}{cccc}
& shellable & & 2-shellable\\
1-dim. &
\begin{tikzpicture}
\coordinate (a) at (0,0);\fill (0,0) circle (1pt);
\coordinate (b) at (0.5,0);\fill (0.5,0) circle (1pt);
\coordinate (c) at (1,0);\fill (1,0) circle (1pt);
\coordinate (d) at (1.5,0);\fill (1.5,0) circle (1pt);
\draw[black] (a) -- (b) -- (c) -- (d) ;
\end{tikzpicture}
& &
\begin{tikzpicture}
\coordinate (a) at (0,0);\fill (0,0) circle (1pt);
\coordinate (b) at (0.5,0);\fill (0.5,0) circle (1pt);
\coordinate (c) at (1,0);\fill (1,0) circle (1pt);
\coordinate (e) at (0.25,0.43);\fill (0.25,0.43) circle (1pt);
\coordinate (f) at (0.75,0.43);\fill (0.75,0.43) circle (1pt);
\coordinate (g) at (1.25,0.43);\fill (1.25,0.43) circle (1pt);
\draw[black] (a) -- (e);
\draw[black] (b) -- (f);
\draw[black] (c) -- (g);
\end{tikzpicture}\\ \\ \\
2-dim. &
\begin{tikzpicture}
\coordinate (a) at (0,0);\fill (0,0) circle (1pt);
\coordinate (b) at (0.5,0);\fill (0.5,0) circle (1pt);
\coordinate (c) at (1,0);\fill (1,0) circle (1pt);
\coordinate (d) at (0.25,0.43);\fill (0.25,0.43) circle (1pt);
\coordinate (e) at (0.75,0.43);\fill (0.75,0.43) circle (1pt);
\draw[black] (a) -- (b) -- (c);
\draw[black] (a) -- (d) -- (b);
\draw[black] (d) -- (e) -- (b);
\draw[black] (e) -- (c);

\colorlet{triangle1}{black!20!white}
\begin{pgfonlayer}{background}
\fill[triangle1!80] (a) -- (b) -- (d) -- cycle;
\fill[triangle1!80] (b) -- (d) -- (e) -- cycle;
\fill[triangle1!80] (b) -- (c) -- (e) -- cycle;
\end{pgfonlayer}
\end{tikzpicture}\quad
\begin{tikzpicture}
\coordinate (a) at (0,0);\fill (0,0) circle (1pt);
\coordinate (b) at (0.8,0);\fill (0.8,0) circle (1pt);
\coordinate (c) at (0.4,0.3);\fill (0.4,0.3) circle (1pt);
\coordinate (d) at (0.4,0.8);\fill (0.4,0.8) circle (1pt);
\draw[black] (a) -- (b) -- (c) -- (a);
\draw[black] (a) -- (c) -- (d) -- (a);
\draw[black] (b) -- (c) -- (d) -- (b);

\colorlet{triangle1}{black!20!white}
\begin{pgfonlayer}{background}
\fill[triangle1!80] (a) -- (b) -- (c) -- cycle;
\fill[triangle1!80] (a) -- (c) -- (d) -- cycle;
\fill[triangle1!80] (b) -- (c) -- (d) -- cycle;
\end{pgfonlayer}
\end{tikzpicture}
& &
\begin{tikzpicture}
\coordinate (a) at (0,0);\fill (0,0) circle (1pt);
\coordinate (b) at (0.5,0);\fill (0.5,0) circle (1pt);
\coordinate (c) at (1,0);\fill (1,0) circle (1pt);
\coordinate (d) at (1.5,0);\fill (1.5,0) circle (1pt);
\coordinate (e) at (0.25,0.43);\fill (0.25,0.43) circle (1pt);
\coordinate (f) at (0.75,0.43);\fill (0.75,0.43) circle (1pt);
\coordinate (g) at (1.25,0.43);\fill (1.25,0.43) circle (1pt);
\draw[black] (a) -- (b) -- (c) -- (d);
\draw[black] (a) -- (e) -- (b);
\draw[black] (b) -- (f) -- (c);
\draw[black] (c) -- (g) -- (d);

\colorlet{triangle1}{black!20!white}
\begin{pgfonlayer}{background}
\fill[triangle1!80] (a) -- (b) -- (e) -- cycle;
\fill[triangle1!80] (b) -- (c) -- (f) -- cycle;
\fill[triangle1!80] (c) -- (d) -- (g) -- cycle;
\end{pgfonlayer}
\end{tikzpicture}\quad
\begin{tikzpicture}
\coordinate (a) at (0.25,0);\fill (0.25,0) circle (1pt);
\coordinate (b) at (0.75,0);\fill (0.75,0) circle (1pt);
\coordinate (c) at (0.5,0.43);\fill (0.5,0.43) circle (1pt);
\coordinate (d) at (0,0.43);\fill (0,0.43) circle (1pt);
\coordinate (e) at (0.25,0.86);\fill (0.25,0.86) circle (1pt);
\coordinate (f) at (0.75,0.86);\fill (0.75,0.86) circle (1pt);
\coordinate (g) at (1,0.43);\fill (1,0.43) circle (1pt);
\draw[black] (a) -- (b) -- (c) -- (a);
\draw[black] (c) -- (d) -- (e) -- (c);
\draw[black] (c) -- (f) -- (g) -- (c);

\colorlet{triangle1}{black!20!white}
\begin{pgfonlayer}{background}
\fill[triangle1!80] (a) -- (b) -- (c) -- cycle;
\fill[triangle1!80] (c) -- (d) -- (e) -- cycle;
\fill[triangle1!80] (c) -- (f) -- (g) -- cycle;
\end{pgfonlayer}
\end{tikzpicture}\\ \\ \\
3-dim. &
\begin{tikzpicture}
\coordinate (a) at (0,0);\fill (0,0) circle (1pt);
\coordinate (b) at (1,0);\fill (1,0) circle (1pt);
\coordinate (c) at (0.5,0.86);\fill (0.5,0.86) circle (1pt);
\coordinate (d) at (1.5,1);\fill (1.5,1) circle (1pt);
\coordinate (e) at (-.5,1);\fill (-0.5,1) circle (1pt);
\coordinate (f) at (0.5,-0.5);\fill (0.5,-0.5) circle (1pt);
\draw[dashed] (a) -- (b) -- (c) -- (a);
\draw[black] (c) -- (d) -- (a) -- (d) -- (b);
\draw[black] (c) -- (e) -- (a) -- (e) -- (b);
\draw[black] (c) -- (f) -- (a) -- (f) -- (b);
\colorlet{triangle1}{black!20!white}
\begin{pgfonlayer}{background}
\fill[triangle1!80] (a) -- (b) -- (c) -- cycle;
\fill[triangle1!80] (a) -- (c) -- (d) -- cycle;
\fill[triangle1!80] (b) -- (c) -- (d) -- cycle;
\fill[triangle1!80] (e) -- (a) -- (c) -- cycle;
\fill[triangle1!80] (f) -- (a) -- (b) -- cycle;
\end{pgfonlayer}
\end{tikzpicture}
\begin{tikzpicture}
\coordinate (a) at (0,0);\fill (0,0) circle (1pt);
\coordinate (b) at (1,0);\fill (1,0) circle (1pt);
\coordinate (c) at (0.5,0.86);\fill (0.5,0.86) circle (1pt);
\coordinate (d) at (1.5,1);\fill (1.5,1) circle (1pt);
\coordinate (e) at (-.5,1);\fill (-0.5,1) circle (1pt);
\coordinate (f) at (2,0.25);\fill (2,0.25) circle (1pt);
\draw[dashed] (b) -- (c);
\draw[dashed] (f) -- (c);
\draw[dashed] (e) -- (b) -- (d);
\draw[black] (b) -- (a) -- (c) -- (d) -- (a) -- (d) -- (b);
\draw[black] (c) -- (e) -- (a) -- (e);
\draw[black] (f) -- (d) -- (f) -- (b);
\colorlet{triangle1}{black!20!white}
\begin{pgfonlayer}{background}
\fill[triangle1!80] (a) -- (b) -- (c) -- cycle;
\fill[triangle1!80] (a) -- (c) -- (d) -- cycle;
\fill[triangle1!80] (b) -- (c) -- (d) -- cycle;
\fill[triangle1!80] (e) -- (a) -- (c) -- cycle;
\fill[triangle1!80] (f) -- (c) -- (d) -- cycle;
\fill[triangle1!80] (f) -- (b) -- (c) -- cycle;
\end{pgfonlayer}
\end{tikzpicture}
 & &
\begin{tikzpicture}
\coordinate (a) at (0,0.25);\fill (0,0.25) circle (1pt);
\coordinate (b) at (1,0.25);\fill (1,0.25) circle (1pt);
\coordinate (c) at (0.5,0.86);\fill (0.5,0.86) circle (1pt);
\coordinate (d) at (1.5,1);\fill (1.5,1) circle (1pt);
\coordinate (e) at (-.25,1);\fill (-0.25,1) circle (1pt);
\coordinate (f) at (1.5,0.25);\fill (1.5,0.25) circle (1pt);
\coordinate (g) at (0.25,-0.25);\fill (0.25,-0.25) circle (1pt);
\coordinate (h) at (0.75,-0.5);\fill (0.75,-0.5) circle (1pt);
\draw[dashed] (f) -- (c);
\draw[dashed] (g) -- (b);
\draw[dashed] (e) -- (b) -- (d);
\draw[black] (b) -- (a) -- (c) -- (d) -- (d) -- (b) -- (c);
\draw[black] (c) -- (e) -- (a) -- (e);
\draw[black] (f) -- (d) -- (f) -- (b);
\draw[black] (c) -- (g) -- (h) -- (c) -- (b) -- (h);
\colorlet{triangle1}{black!20!white}
\begin{pgfonlayer}{background}
\fill[triangle1!80] (a) -- (b) -- (c) -- cycle;
\fill[triangle1!80] (a) -- (c) -- (d) -- cycle;
\fill[triangle1!80] (b) -- (c) -- (d) -- cycle;
\fill[triangle1!80] (e) -- (a) -- (c) -- cycle;
\fill[triangle1!80] (f) -- (c) -- (d) -- cycle;
\fill[triangle1!80] (f) -- (b) -- (c) -- cycle;
\fill[triangle1!80] (g) -- (b) -- (c) -- cycle;
\fill[triangle1!80] (g) -- (h) -- (b) -- cycle;
\end{pgfonlayer}
\end{tikzpicture}\
\begin{tikzpicture}
\coordinate (a) at (-0.5,0);\fill (-0.5,0) circle (1pt);
\coordinate (b) at (0,0);\fill (0,0) circle (1pt);
\coordinate (c) at (0,0.75);\fill (0,0.75) circle (1pt);
\coordinate (d) at (-0.75,1);\fill (-0.75,1) circle (1pt);
\coordinate (e) at (0.5,0.75);\fill (0.5,0.75) circle (1pt);
\coordinate (f) at (0.75,-0.25);\fill (0.75,-0.25) circle (1pt);
\coordinate (g) at (1.25,0);\fill (1.25,0) circle (1pt);
\coordinate (h) at (1.25,0.5);\fill (1.25,0.5) circle (1pt);
\draw[dashed] (b) -- (e);
\draw[dashed] (a) -- (c);
\draw[dashed] (f) -- (h);
\draw[black] (b) -- (d) -- (a) -- (b) -- (c) -- (d);
\draw[black] (c) -- (e) -- (f) -- (b) -- (c) -- (f) -- (g);
\draw[black] (f) -- (e) -- (h) -- (g) -- (e);
\colorlet{triangle1}{black!20!white}
\begin{pgfonlayer}{background}
\fill[triangle1!80] (a) -- (b) -- (c) -- cycle;
\fill[triangle1!80] (b) -- (c) -- (e) -- cycle;
\fill[triangle1!80] (b) -- (e) -- (f) -- cycle;
\fill[triangle1!80] (a) -- (c) -- (d) -- cycle;
\fill[triangle1!80] (f) -- (h) -- (e) -- cycle;
\fill[triangle1!80] (f) -- (h) -- (g) -- cycle;
\end{pgfonlayer}
\end{tikzpicture}\
\begin{tikzpicture}
\coordinate (a) at (-1,-0.25);\fill (-1,-0.25) circle (1pt);
\coordinate (b) at (0,0.2);\fill (0,0.2) circle (1pt);
\coordinate (c) at (0,0.75);\fill (0,0.75) circle (1pt);
\coordinate (d) at (-0.5,0.8);\fill (-0.5,0.8) circle (1pt);
\coordinate (e) at (0.5,0.75);\fill (0.5,0.75) circle (1pt);
\coordinate (f) at (0.75,-0.25);\fill (0.75,-0.25) circle (1pt);
\coordinate (g) at (-0.1,1.25);\fill (-0.1,1.25) circle (1pt);
\coordinate (h) at (-1,2);\fill (-1,2) circle (1pt);
\draw[dashed] (d) -- (b) -- (e);
\draw[black] (a) -- (c);
\draw[dashed] (b) -- (d) -- (g);
\draw[black] (d) -- (h) -- (g);
\draw[black] (d) -- (a) -- (b) -- (c) -- (d);
\draw[black] (h) -- (c) -- (e) -- (f) -- (b) -- (c) -- (f);
\draw[black] (c) -- (g);
\colorlet{triangle1}{black!20!white}
\begin{pgfonlayer}{background}
\fill[triangle1!80] (a) -- (b) -- (c) -- cycle;
\fill[triangle1!80] (b) -- (c) -- (e) -- cycle;
\fill[triangle1!80] (b) -- (e) -- (f) -- cycle;
\fill[triangle1!80] (a) -- (c) -- (d) -- cycle;
\fill[triangle1!80] (d) -- (h) -- (g) -- cycle;
\fill[triangle1!80] (d) -- (h) -- (c) -- cycle;
\fill[triangle1!80] (c) -- (h) -- (g) -- cycle;
\end{pgfonlayer}
\end{tikzpicture}
\end{array}$$
\caption{}\label{shell exams}
\end{figure}
\end{exam}

\begin{thm}\label{link}
Let $\Del$ be a $k$-shellable complex and $\sigma$ a face of $\Del$. Then $\lk_\Del(\sigma)$ is again $k$-shellable.
\end{thm}
\begin{proof}
Since $\Del$ is $k$-shellable, so there exists an $k$-shelling order $F_1,\ldots,F_r$ of facets of $\Del$. Let $F_{i_1},\ldots,F_{i_t}$ where $i_1<\ldots<i_t$ be all of facets which contain $\sigma$. We claim that $F_{i_1}\backslash \sigma,\ldots,F_{i_t}\backslash \sigma$ is a $k$-shelling order of $\lk_\Del(\sigma)$. To this end we want to show that the condition (c) of Proposition \ref{k-shell} holds.

Set $G_j=F_{i_j}\backslash \sigma$. Consider $l,m$ with $1\leq l< m\leq t$. By $k$-shellability of $\Del$, there are
$x_1,\ldots,x_k\in F_{i_m}\backslash F_{i_l}=(F_{i_m}\backslash \sigma)\backslash (F_{i_l}\backslash \sigma)=G_m\backslash G_l$ such that for some $s< i_m$ we have $\{x_1,\ldots,x_k\}=F_{i_m}\backslash F_s$. It follows from $\sigma\subset F_{i_m}$ and $F_{i_m}\backslash F_s=\{x_1,\ldots,x_k\}$ that $\sigma\subset F_s$. This implies that $F_s$ is among the list $F_{i_1},\ldots,F_{i_m}$. Let $F_{i_{m'}}=F_s$. Hence $G_m\backslash G_{m'}=\{x_1,\ldots,x_k\}$ and the assertion is completed.
\end{proof}

For the simplicial complexes
$\Del_1$ and $\Del_2$ defined on disjoint vertex sets, the \emph{join}
of $\Del_1$ and $\Del_2$ is $\Del_1\cdot\Del_2=\{\sigma\cup\tau:\
\sigma\in\Del_1,\tau\in\Del_2\}$.

\begin{thm}
The simplicial complexes $\Del_1$ and $\Del_2$ are $k$-shellable if and only if $\Del_1\cdot\Del_2$ is $k$-shellable.
\end{thm}
\begin{proof}
Let $\Del_1$ and $\Del_2$ be $k$-shellable. Let $F_1,\ldots,F_r$ and $G_1,\ldots,G_s$ be, respectively, the $k$-shelling orders of $\Del_1$ and $\Del_2$. We claim that 
$$F_1\cup G_1,F_1\cup G_2,\ldots,F_1\cup G_s,\ldots,F_r\cup G_1,F_r\cup G_2,\ldots,F_r\cup G_s$$
is a $k$-shelling order of $\Del_1\cdot\Del_2$.

Let $F_i\cup G_j$ be a facet of $\Del_1\cdot\Del_2$ which comes after $F_p\cup G_q$ in the above order. We have some cases:

Let $p<i$. Since $\Del_1$ is $k$-shellable, there exist $u_1,\ldots,u_k\in F_i\backslash F_p$ and some $l<i$ such that $F_i\backslash F_l=\{u_1,\ldots,u_k\}$. It follows that $u_1,\ldots,u_k\in (F_i\cup G_j)\backslash (F_p\cup G_q)$ and $(F_i\cup G_j)\backslash (F_l\cup G_j)=\{u_1,\ldots,u_k\}$.

Let $p=i$ and $q<j$. Since $\Del_2$ is $k$-shellable, there exist $v_1,\ldots,v_k\in G_j\backslash G_q$ and some $m<j$ such that $G_j\backslash G_m=\{v_1,\ldots,v_k\}$. Therefore we obtain $v_1,\ldots,v_k\in (F_i\cup G_j)\backslash (F_p\cup G_q)$ and $(F_i\cup G_j)\backslash (F_i\cup G_m)=\{v_1,\ldots,v_k\}$.

Conversely, suppose that $\Del_1\cdot\Del_2$ is $k$-shellable with the $k$-shelling order $F_{i_1}\cup G_{j_1},\ldots,F_{i_t}\cup G_{j_t}$. Let $F_{s_1},\ldots,F_{s_r}$ be the ordering obtained from $F_{i_1}\cup G_{j_1},\ldots,F_{i_t}\cup G_{j_t}$ after removing the repeated facets beginning on the left-hand. Then it is easy to check that $F_{s_1},\ldots,F_{s_r}$ is a $k$-shelling order of $\Del_1$. In a similar way, it is shown that $\Del_2$ is $k$-shellable.
\end{proof}

The following theorem, relates the expansion of a shellable complex to a $k$-shellable complex.

\begin{thm}\label{sh and k-sh}
Let $\Del$ be a simplicial complex and $k\in\NN$. Then $\Del$ is shellable if and only if $\Del^\k$ is $k$-shellable.
\end{thm}
\begin{proof}
Let $\Del=\langle F_1,\ldots,F_r\rangle$ and let $\Del_j=\langle F_j\rangle\cap\langle F_1,\ldots,F_{j-1}\rangle$ for $j=2,\ldots,r$. Fix an integer $j$. If
$\Del_j=\langle F_j\backslash x_{i_1},\ldots,F_j\backslash x_{i_t}\rangle$, then
$$\begin{array}{rl}
    \Del^\k_j= & \langle F^\k_j\rangle\cap\langle F^\k_1,\ldots,F^\k_{j-1}\rangle \\
    =& \langle F^\k_j\backslash \{x_{i_11},\ldots,x_{i_1k}\},\ldots,F^\k_j\backslash
\{x_{i_t1},\ldots,x_{i_tk}\}\rangle.
  \end{array}
$$
Now by the Definition \ref{Def k-shell}, if $F_1,\ldots,F_r$ is a shelling order of $\Del$ then $F^\k_1,\ldots,F^\k_r$ is a $k$-shelling order
of $\Del^\k$.

Conversely, suppose that $F^\k_1,\ldots,F^\k_r$ is a $k$-shelling order of $\Del^\k$ and set $\Del^\k_j=\langle F^\k_j\rangle\cap\langle F^\k_1,\ldots,F^\k_{j-1}\rangle$ for $j=2,\ldots,r$. Fix an index $j$. Hence
$\Del^\k_j=\langle F^\k_j\backslash \sigma_1,\ldots,F^\k_j\backslash\sigma_t\rangle$
with $|\sigma_i|=k$ for all $i$. By Proposition \ref{k-shell}(b), $\sigma_l\cap\sigma_m=\emptyset$ for all $l\neq m$. We claim that for every $i$, $\sigma_i$ is the expansion of a singleton set. Suppose, on the contrary, that for some $\sigma_s$ we have $x_{i_1l},x_{i_2m}\in \sigma_s$ with $i_1\neq i_2$ and let $F^\k_j\cap F^\k_{s'}=F^\k_j\backslash \sigma_s$ for some $s'$. It follows from $|\sigma_s|=k$ that $x_{i_1l'}\not\in\sigma_s$ for some $l'$ with $1\leq l'\leq k$. In particular, we conclude that $x_{i_1l}\not\in F^\k_{s'}$ but $x_{i_1l'}\in F^\k_{s'}$. This is a contradiction, because $F^\k_{s'}$ is the expansion of $F_{s'}$.

Therefore we conclude that for all $j=2,\ldots,r$, the complex $\Del^\k_j$ is in the form
$$\Del^\k_j=\langle F^\k_j\backslash\{x_{i_1}\}^\k,\ldots,F^\k_j\backslash\{x_{i_t}\}^\k\rangle.$$
Finally, for all $j$, $\Del_j=\langle F_j\rangle\cap\langle F_1,\ldots,F_{j-1}\rangle$ will be in the form $\Del_j=\langle F_j\backslash x_{i_1},\ldots,F_j\backslash x_{i_t}\rangle$. This implies that $\langle F_j\rangle\cap\langle F_1,\ldots,F_{j-1}\rangle$ is pure of dimension $\dim(F_j)-1$ for all $j\geq 2$, as desired.
\end{proof}

Let $R$ be a Noetherian ring and $M$ be a finitely generated multigraded $R$-module. We call
$$\F:0=M_0\subset M_1\subset \ldots\subset M_{r-1}\subset M_r=M$$
a \emph{multigraded finite filtration} of submodules of $M$ if there exist the positive integers $a_1,\ldots,a_r$ such that $M_i/M_{i-1}\cong \prod^{a_i}_{j=1}R/P_i(-\a_{ij})$ for some $P_i\in\Supp(M)$. $\F$ is called a \emph{multigraded prime filtration} if $a_1=\ldots=a_r=1$. It is well known that every finitely generated multigraded $R$-module $M$ has a multigraded prime filtration (see for example \cite[Theorem 6.4]{Ma}). In the following we present a multigraded finite filtration for the face ring of a $k$-shellable simplicial complex which we need in Section 3.

For $F\subset X$. We set $F^c=X\backslash F$ and $P_F=(x_i:x_i\in F)$.

\begin{thm}\label{k-sh and k-cl}
Let $\Del$ be a simplicial complex and $k$ a positive integer. If $F_1,\ldots,F_r$ is a $k$-shelling order of $\Del$
then there exists a filtration $0=M_0\subset M_1\subset\ldots\subset M_r=S/I_\Del$ with
$$M_i=\bigcap^{r-i}_{l=1}P_{F^c_l}\quad and\quad M_i/M_{i-1}\cong\prod^{k^{a_i}}_{j=1} S/P_{F^c_{r-i+1}}(-\a_{ij}),$$
for all $i=1,\ldots,r$. Here $a_i=|\a_{ij}|$ for all $j=1,\ldots,k^{a_i}$.
\end{thm}
\begin{proof}
We set $a_1=0$ and for each $i>2$ we denote by $a_i$ the number of facets of $\langle F_i\rangle\cap\langle F_1,\ldots,F_{i-1}\rangle$. If $F_1,\ldots,F_r$ is a $k$-shelling of $\Del$, then for $i=2,\ldots,r$ we have
\begin{equation}\label{1}
\bigcap^{i-1}_{j=1}P_{F^c_j}+P_{F^c_i}=P_{F^c_i}+P_{\sigma_{i1}}\ldots P_{\sigma_{ia_i}}
\end{equation}
where $\sigma_{il}=F_i\backslash F_{i_l}$ and $|\sigma_{il}|=k$ for $l=1,\ldots,a_i$. Actually, $ F_{i_l}\cap F_i$'s
are all of facets of $\langle F_i\rangle\cap\langle F_1,\ldots,F_{i-1}\rangle$. Since
$\sigma_{il}\cap\sigma_{il'}=\emptyset$ for $1\leq l<l'\leq a_i$, one can suppose that
$P_{\sigma_{i1}}\ldots P_{\sigma_{ia_i}}=(f_{ij}:j=1,\ldots,k^{a_i})$. Set $\a_{ij}=\deg(f_{ij})$ and it is clear that for
all $j=1,\ldots,k^{a_i}$, $a_i=|\a_{ij}|$. We have the following isomorphisms:
$$\begin{array}{ll}
    (\bigcap^{i-1}_{j=1}P_{F^c_j})/(\bigcap^{i}_{j=1}P_{F^c_j}) & \cong\bigcap^{i-1}_{j=1}P_{F^c_j}+P_{F^c_i}/P_{F^c_i}
    \\
    &\cong P_{\sigma_{i1}}\ldots P_{\sigma_{ia_i}}+P_{F^c_i}/P_{F^c_i} \\
     & \cong P_{\sigma_{i1}}\ldots P_{\sigma_{ia_i}}/(P_{\sigma_{i1}}\ldots P_{\sigma_{ia_i}}\cap P_{F^c_i})\\
     & \cong P_{\sigma_{i1}}\ldots P_{\sigma_{ia_i}}/P_{\sigma_{i1}}\ldots P_{\sigma_{ia_i}}P_{F^c_i}
  \end{array}
$$
where $a_i=|\a_{ij}|$ for $j=1,\ldots,k^{a_i}$. Now it is easy to check that the homomorphism
$$\begin{array}{rcc}
    \theta:\prod^{k^{a_i}}S&\rightarrow & P_{\sigma_{i1}}\ldots
    P_{\sigma_{ia_i}}/P_{\sigma_{i1}}\ldots P_{\sigma_{ia_i}}P_{F^c_i} \\
     (r_1,\ldots,r_{k^{a_i}}) & \mapsto & \sum^{k^{a_i}}_jr_jf_{ij}+P_{\sigma_{i1}}\ldots P_{\sigma_{ia_i}}P_{F^c_i}
  \end{array}
$$
is an epimorphism. In particular, it follows that
$$P_{\sigma_{i1}}\ldots P_{\sigma_{ia_i}}/P_{\sigma_{i1}}\ldots P_{\sigma_{ia_i}}P_{F^c_i}\cong\prod^{k^{a_i}}_{j=1} S/P_{F^c_i}(-\a_{ij}).$$
This completes the proof.
\end{proof}

\begin{rem}
In view of Theorem \ref{k-sh and k-cl}, let $\Del=\langle F_1,\ldots,F_r\rangle$ be a shellable simplicial complex and
$\Del_j=\langle F_1,\ldots,F_j\rangle$. Then we have the prime filtration
$$(0)=I_\Del\subset I_{\Del_{r-1}}\subset\ldots\subset I_{\Del_1}\subset K[\Del]$$
for $K[\Del]$. In particular, it follows the following filtration for $K[\Del^\k]$:
$$(0)=I_{\Del^\k}\subset I_{\Del^\k_{r-1}}\subset\ldots\subset I_{\Del^\k_1}\subset K[\Del^\k].$$
In other words, Theorem \ref{k-sh and k-cl} gives a filtration for the face ring of the expansion of a shellable
simplicial complex with respect to $k$.
\end{rem}

\begin{rem}
The filtration described in Theorem \ref{k-sh and k-cl} in the case that $k>1$ is not a prime filtration, i.e. the quotient of any two consecutive modules of the filtration is not cyclic. Consider the same notations of Theorem \ref{k-sh and k-cl}, we have the following prime filtration for $K[\Del]$ when $\Del$ has a $k$-shelling order:
\begin{flushleft}
$\F: 0=\bigcap^r_{j=1}P_{F^c_j}\subset\ldots\subset$
\end{flushleft}
\begin{center}
$\bigcap^{i}_{j=1}P_{F^c_j}\subset\ldots\subset \sum^{k^{a_i}-1}_{j=1}(f_j)+\bigcap^{i}_{j=1}P_{F^c_j}\subset
\sum^{k^{a_i}}_{j=1}(f_j)+\bigcap^{i}_{j=1}P_{F^c_j}=\bigcap^{i-1}_{j=1}P_{F^c_j}$
\end{center}
\begin{flushright}
$\subset\ldots\subset K[\Del]$
\end{flushright}
where $(f_1,\ldots,f_{j-1}):(f_j)$ is generated by linear forms for all $j=2,\ldots,k^{a_i}$ and all $i=1,\ldots,r$.

For all $2\leq j\leq k^{a_i}$, suppose that $(f_1,\ldots,f_{j-1}):(f_j)=P_{Q_j}$. Set $P_{Q_1}=(0)$. We have
$$\begin{array}{rl}
    \sum^{j}_{t=1}(f_t)+\bigcap^{i}_{t=1}P_{F^c_t}/ \sum^{j-1}_{t=1}(f_t)+\bigcap^{i}_{t=1}P_{F^c_t}\cong&
    (f_j)/(f_j)\cap\big((f_1,\ldots,f_{j-1})+\bigcap^{i}_{t=1}P_{F^c_t}\big) \\
     \cong& (f_j)/f_jP_{L_{ij}}.
  \end{array}
$$
where $L_{ij}=F^c_i\cup Q_j$, for $i=1,\ldots,r$ and $j=1,\ldots,k^{a_i}$.
Therefore the set of prime ideals which defines the cyclic quotients of $\F$ is $\Supp(\F)=\{P_{L_{ij}}:i=1,\ldots,r\ \mbox{and}\ j=1,\ldots,k^{a_i}\}$.
\end{rem}

\section{The Stanley conjecture}

Consider a field $K$, and let $R$ be a finitely generated $\NN^n$-graded $K$-algebra, and let $M$ be a finitely generated $\ZZ^n$-graded $R$-module. Stanley \cite{St} conjectured that, in this case, there exist finitely many subalgebras $A_1,\ldots,A_r$ of $R$,
each generated by algebraically independent $\NN^n$-homogeneous elements of $R$, and there exist $\ZZ^n$-homogeneous elements
$u_1,\ldots,u_r$ of $M$, such that $M=\bigoplus^r_{i=1}u_iA_i$, where $\dim(A_i)\geq \depth(M)$ for all $i$ and where
$u_iA_i$ is a free $A_i$-module of rank one.

Consider a finitely generated $\ZZ^n$-graded $S$-module $M$, a subset $Z$ of
$\{x_1,\ldots,x_n\}$, and a homogeneous element $u\in M$. The $K$-subspace $uK[Z]$ of $M$ is called a \emph{Stanley space of dimension
$|Z|$} if it is a free $K[Z]$-module of rank 1, i.e., the elements of the form $uv$, where $v$ is a monomial in $K[Z]$, form a $K$-basis of $uK[Z]$. A \emph{Stanley decomposition} of $M$ is a decomposition $\D$ of $M$ into a finite direct sum of Stanley spaces. The \emph{Stanley depth}
of $\D$, denoted $\sdepth(\D)$, is the minimal dimension of a Stanley space in a decomposition $\D$. We set
$$\sdepth(M )=\max\{\sdepth(\D):\D\ \mbox{is a Stanley decomposition of}\ M\},$$
and we call this number the \emph{Stanley depth of $M$}. The Stanley conjecture says that $\sdepth(M)\geq\depth(M)$ always holds.

The following lemma is needed in the proof of the main theorem of this section.

\begin{lem}\label{product}
Let $F_1,\ldots,F_s\subset X$ with $F_i\cap F_j=\emptyset$ for all $i\neq j$ and $|F_i|=k>1$. Let
$f_1,\ldots,f_{k^s}$ be a sequence of minimal generators of $P_{F_1}\ldots P_{F_s}$ ordered with respect to
lexicographical ordering $x_1>x_2>\ldots>x_n$. Suppose that $n_i$ is the minimal number of homogeneous generators of
$(f_1,\ldots,f_{i-1}):(f_i)$ for $i=2,\ldots,k^s$. Then
$$\max\{n_i:i=2,\ldots,k^s\}=n_{k^s}=(k-1)s.$$
Moreover, for all $i$, the colon ideal $(f_1,\ldots,f_{i-1}):(f_i)$ is generated by linear forms.
\end{lem}
\begin{proof}
By \cite[Corollary 1.5.]{CoHe}, $P_{F_1}\ldots P_{F_s}$ has linear quotients. To show equality, we use induction on $s$. If $s=1$, the assertion is clear. Assume that $s>1$. Let $P_{F_1}\ldots P_{F_{s-1}}=(f_1,\ldots,f_{k^{s-1}})$ and $F_s=\{x_1,\ldots,x_k\}$. Then
$$P_{F_1}\ldots P_{F_s}=x_1(f_1,\ldots,f_{k^{s-1}})+\ldots +x_k(f_1,\ldots,f_{k^{s-1}}).$$
Moreover, $$(x_1f_1,\ldots,x_1f_{k^{s-1}},\ldots,x_kf_1,\ldots,x_kf_{k^{s-2}}):x_kf_{k^{s-1}}=(x_1,\ldots,x_{k-1})+(f_1,\ldots,f_{k^{s-2}}):f_{k^{s-1}}.$$
Now, by the induction hypothesis, we have $n_{k^s}=(k-1)+(k-1)(s-1)=(k-1)s$, as desired.
\end{proof}

Let
$$\F:0=M_0\subset M_1\subset \ldots\subset M_r=M$$
be a prime filtration of $M$ with $M_i/M_{i-1}\cong (S/P_i)(-\a_i)$. Then this filtration
decomposes $M$ as a multigraded $K$-vector space, that is, we have $M=\bigoplus^r_{i=1}u_iK[Z_i]$ and this is a Stanley
decomposition of $M$ where $u_i\in M_i$ is a homogeneous element of degree $\a_i$ and $Z_i=\{x_j:x_j\not\in P_i\}$.

Now suppose that $\Del=\langle G_1,\ldots,G_r\rangle$ is a pure shellable simplicial complex on $X$ and $k$ is a
positive integer. By Theorem \ref{sh and k-sh}, $\Del^\k=\langle G^\k_1,\ldots,G^\k_r\rangle$ is pure $k$-shellable.
For all $i$, set $F_i=G^\k_i$. Consider the prime filtration $\F$ of $K[\Del^\k]$ described in Section 2. Then we have
the Stanley decomposition $$K[\Del^\k]=\bigoplus^r_{i=1}\bigoplus^{k^{a_i}}_{j=1}u_{ij}K[Z_{ij}]$$
where $Z_{ij}=\{x_l:l\not\in L_{ij}\}$, $\deg(u_{ij})=\a_{ij}$ and $|\a_{ij}|=a_i$ for all $i,j$. We claim that for all
$i,j$, $|Z_{ij}|\geq\depth(K[\Del^\k])$.

By Corollaries 4.1 and 2.1 of \cite{Ra}, we have $\depth(K[\Del^\k])=\dim(K[\Del])=|G_i|$. On the other hand,
$|G_i|\geq a_i$. Now by combining all of these results with Lemma \ref{product}, we have
$$\begin{array}{rl}
   |Z_{ik^{a_i}}| & =kn-(|Q_{k^{a_i}}|+|F^c_i|) \\
     & =kn-((k-1)a_i+k\h(P_{G^c_i})) \\
     & =(k-1)(n-\h(P_{G^c_i})-a_i)+n-\h(P_{G^c_i}) \\
     & = (k-1)(|G_i|-a_i)+|G_i|\\
     & \geq |G_i|=\dim(K[\Del]).
  \end{array}
$$
Thus we have shown the main result of this section:

\begin{thm}\label{k-sh stanley}
The expansion of the face ring of a pure shellable simplicial complex with respect to $k>0$ satisfies the Stanley
conjecture. In particular, the face ring of a pure $k$-shellable complex satisfies the Stanley conjecture.
\end{thm}

\section{$k$-shellable graphs}

Let $G$ be a simple graph and let $\Del_G$ the independence complex of $G$. We say that $G$ is $k$-shellable if $\Del_G$ has this property. The purpose of this section is to characterize $k$-shellable graphs.

Following Schrijver \cite{Sh}, the \emph{duplication} of a vertex $x_i$ of a graph $G$ means extending its vertex set $X$ by a new vertex $x_{i'}$ and replacing $E(G)$ by $$E(G)\cup\{x_{i'}x_j:x_ix_j\in E(G)\}.$$

In other words, if $V(G)=\{x_1,\ldots,x_n\}$ then the graph $G'$ obtained from $G$ by duplicating $k_i-1$ times the vertex $x_i$ has the vertex set
$$V(G')=\{x_{ij}:i=1,\ldots,n\ \mbox{and}\ j=1,\ldots.k_i\}$$
and the edge set
$$E(G')=\{x_{ir}x_{js}:x_ix_j\in E(G),\ r=1,\ldots,k_i\ \mbox{and}\ j=1,\ldots,k_i\}.$$

\begin{exam}
Let $G$ be a simple graph on the vertex set $V(G)=\{x_1,\ldots,x_5\}$ and $E(G)=\{x_1x_3,x_1x_4,x_2x_4,x_2x_5,x_3x_5,x_4x_5\}$. Let $G'$ be obtained from $G$ by duplicating $1$ times the vertices $x_1$ and $x_4$ and $0$ times the other vertices. Then $G$ and $G'$ are in the form
$$\begin{array}{ccccc}
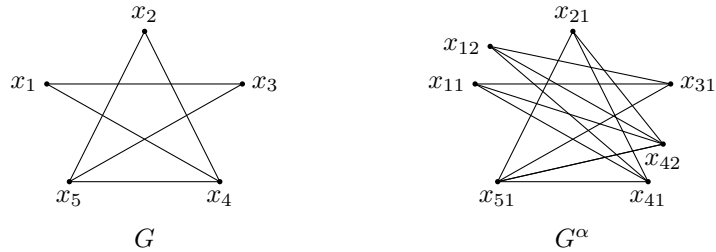
\begin{tikzpicture}
\coordinate (x1) at (-1.3,1.3);\fill (-1.3,1.3) circle (1pt);
\node[left] at (-1.3,1.3)    {$x_1$};
\coordinate (x2) at (0,2);\fill (0,2) circle (1pt);
\node[above] at (0,2)    {$x_2$};
\coordinate (x3) at (1.3,1.3);\fill (1.3,1.3) circle (1pt);
\node[right] at (1.3,1.3)    {$x_3$};
\coordinate (x4) at (1,0);\fill (1,0) circle (1pt);
\node[below] at (1,0)    {$x_4$};1
\coordinate (x5) at (-1,0);\fill (-1,0) circle (1pt);
\node[below] at (-1,0)    {$x_5$};
\node[below] at (0,-0.5)    {$G$};
\draw (x1) -- (x3);
\draw (x3) -- (x5);
\draw (x5) -- (x2);
\draw (x2) -- (x4) -- (x1);
\draw (x4) -- (x5);
\end{tikzpicture}
& \ & \ & \ &
\begin{tikzpicture}
\coordinate (x11) at (-1.3,1.3);\fill (-1.3,1.3) circle (1pt);
\node[left] at (-1.3,1.3)    {$x_{11}$};
\coordinate (x21) at (0,2);\fill (0,2) circle (1pt);
\node[above] at (0,2)    {$x_{21}$};
\coordinate (x31) at (1.3,1.3);\fill (1.3,1.3) circle (1pt);
\node[right] at (1.3,1.3)    {$x_{31}$};
\coordinate (x41) at (1,0);\fill (1,0) circle (1pt);
\node[below] at (1,0)    {$x_{41}$};1
\coordinate (x51) at (-1,0);\fill (-1,0) circle (1pt);
\node[below] at (-1,0)    {$x_{51}$};
\coordinate (x12) at (-1.1,1.8);\fill (-1.1,1.8) circle (1pt);
\node[left] at (-1.1,1.8)    {$x_{12}$};
\coordinate (x42) at (1.2,0.5);\fill (1.2,0.5) circle (1pt);
\node[below] at (1.2,0.5)    {$x_{42}$};
\node[below] at (0,-0.5)    {$G^\alpha$};
\draw (x11) -- (x31);
\draw (x31) -- (x51);
\draw (x51) -- (x21);
\draw (x21) -- (x41) -- (x11);
\draw (x41) -- (x51);
\draw (x12) -- (x31);
\draw (x41) -- (x12);
\draw (x21) -- (x42) -- (x11);
\draw (x42) -- (x51);
\draw (x42) -- (x12);
\draw (x42) -- (x51);
\end{tikzpicture}
\end{array}$$
Also, the independence complexes of $G$ and $G'$ are, respectively, $\Del_G=\langle x_1x_2,x_1x_5,x_2x_3,x_3x_4\rangle$ and $\Del_{G'}=\langle x_{11}x_{12}x_{21},x_{11}x_{12}x_{51},x_{21}x_{31},x_{31}x_{41}x_{42}\rangle$. Note that $\Del_{G'}=\Del^\alpha_{G}$ where $\alpha=(2,1,1,2,1)$.
\end{exam}

In the following theorem we show that the simple graph obtained from duplicating $k-1$ times any vertex of a shellable graph is $k$-shellable.

\begin{thm}\label{Sh to k-sh}
Let $G$ be a simple graph on $X$ and let $G'$ be a new graph obtained from $G$ by duplicating $k-1$ times any vertex of $G$. Then $G$ is shellable if and only if $G'$ is $k$-shellable.
\end{thm}
\begin{proof}
It suffices to show that $\Del_{G'}=\Del^\k_G$. Then Theorem \ref{sh and k-sh} completes the assertion.

After relabeling of the vertices of $G'$ one can assume that $G'$ is a graph with the vertex set $X^\k=\{x_{ij}:i=1,\ldots,n,\ j=1,\ldots,k\}$ and the edge set $$\{x_{ir}x_{js}:x_ix_j\in G\ \mbox{and}\ 1\leq r,s\leq k\}.$$
Let $F$ be an independent set of $G'$ and let $\bar{F}=\{x_i:x_{ir}\in F\ \mbox{for some}\ r\}$. If $|\bar{F}|=1$ then $\bar{F}$ is an independent set in $G$. So assume that $|\bar{F}|>1$. Suppose, on the contrary, that $x_i,x_j\in \bar{F}$ and $x_ix_j\in G$. By the construction of $G'$, for all $x_{ir}$ and $x_{js}$ of $V(G')$, $x_{ir}x_{js}\in G'$. Therefore $F$ contains an edge $x_{ir}x_{js}$ of $G'$, a contradiction. This implies that $\bar{F}$ is an independent set in $G$. In particular, since $F\subset (\bar{F})^\k$ we have $F\in \Del^\k_G$.

Conversely, suppose $H$ is an independent set of $G$. Choose $x_{ir},x_{js}\in H^\k$. If $x_{ir}x_{js}\in G'$ then $x_ix_j\in G$, which is false since $x_i,x_j\in H$. Therefore $H^\k$ is an independent set of $G'$ and $H^\k\in\Del_{G'}$.
\end{proof}

In the following we want to extend some results from \cite{CaCr,CrEs,VaVi}. Firstly, we present a generalization of the concept of simplicial vertex.

Let $G$ be a simple graph. For $U\subset V(G)$ we define the \emph{induced subgraph} of $G$ on $U$ to be the subgraph $G_U$ on $U$ consisting of those edges $x_ix_j\in E(G)$ with $x_i,x_j\in U$. For $x\in V(G)$, let $N_G(x)$ denote the open neighborhood
of $x$, that is, all of vertices adjacent to $x$. We also denote by $N_G[x]$ the closed neighborhood
of $x$, which is $N_G(x)$ together with $x$ itself, so that $N_G[x]=N_G(x)\cup\{x\}$. Set $N_G(U)=\bigcup_{x\in U}N_G(x)$ and $N_G[U]=\bigcup_{x\in U}N_G[x]$.

Recall from \cite{Di} or \cite{LeBo} that a vertex $x\in V(G)$ is \emph{simplicial} if the induced subgraph $G_{N_G[x]}$ is complete.

The simple graph $G$ is a \emph{complete $r$-partite} graph if there is a partition $V(G)=V_1\cup\ldots\cup V_r$ of the vertex set, such that $uv\in E(G)$ if and only if $u$ and $v$ are in different parts of the partition. If $|V_i|=n_i$, then $G$ is denoted by $K_{n_1.\ldots,n_r}$.

\begin{defn}
The set $S$ of pairwise non-adjacent vertices of $G$ is a \emph{$k$-simplicial set} if $G_{N_G[S]}$ is a $r$-partite complete graph with $k$-element parts $S_1,\ldots,S_r$ having the following property:
\begin{center}
    for every $S_l$ and every two vertices $x_i,x_j\in S_l$, $N_G(x_i)=N_G(x_j)$.
\end{center}
\end{defn}

Note that every $1$-simplicial set is a simplicial vertex.

\begin{lem}\label{ney and lk}
Let $S\subset V(G)$ of pairwise non-adjacent vertices of $G$. Set $G'=G\backslash N_G[S]$ and $G''=G\backslash S$. Then
\begin{center}
    (i) $\Del_{G'}=\lk_{\Del_G}(S)$ and (ii) $\Del_{G''}=\bigcap_{x\in S} \dl_{\Del_G}(x)$.
\end{center}
\end{lem}
\begin{proof}
(i) ``$\subseteq$'': Let $F\in\Del_{G'}$. Since $F\subset V(G)\backslash N_G[S]$ we have $F\in\Del_G$ and $F\cap S=\emptyset$. It remains to show that $F\cup S\in\Del_G$. Let $F\cup S$ contain an edge $x_ix_j\in E(G)$. Then it should be $x_i\in S$ and $x_j\in F$. In particular, since $x_i$ and $x_j$ are adjacent it follows that $x_j\in N_G[S]$. This is impossible because $F\subset V(G)\backslash N_G[S]$. Therefore $F\cup S\in\Del_G$. This implies that $F\in\lk_{\Del_G}(S)$.

``$\supseteq$'': Let $F\in\lk_{\Del_G}(S)$. Then $F\cup S\in\Del_G$ and $F\cap S=\emptyset$. In order to prove that $F\in\Del_{G'}$ it suffices to show that $F\cap N_G(S)=\emptyset$ and no two vertices of $F$ are adjacent in $G'$. If $x\in F\cap N_G(S)$ then $S\cup\{x\}$ contains an edge of $G$ and this contradicts $F\cup S\in\Del_G$. Also, if $x_i,x_j\in F$ and $x_ix_j\in G'$ then $x_ix_j\in G$, which is again a contradiction.

(ii) ``$\subseteq$'': Let $F\in\Del_{G''}$. If for $x_i,x_j\in F$, $x_ix_j$ is an edge of $G$, since $\{x_i,x_j\}\cap S=\emptyset$ we obtain that $x_ix_j\in G''$, which is not true. Hence $F\in\Del_G$. In particular, it follows from $F\subset V(G)\backslash S$ that $F\cap S=\emptyset$ and so $F\in\bigcap_{x\in S} \dl_{\Del_G}(x)$.

``$\supseteq$'': Let $F\in\bigcap_{x\in S} \dl_{\Del_G}(x)$. Then for all $x\in S$ we have $F\in\dl_{\Del_G}(x)$. Thus $F\in\Del_G$ and $F\cap S=\emptyset$. Since $F$ contains no edge of $G$ it follows that $F$ contains no edge of $G''$, either. This implies that $F\in\Del_{G''}$.
\end{proof}

Combining Theorem \ref{link} with Lemma \ref{ney and lk} we obtain the following corollary as an extension of Theorem 2.6. of \cite{VaVi}.

\begin{cor}\label{G to lk}
Let $S\subset V(G)$ be a set of pairwise non-adjacent vertices. If $G$ is $k$-shellable, then $G'=G\backslash N_G[S]$ is $k$-shellable, too.
\end{cor}

\begin{thm}
If $S$ is a $k$-simplicial set of $G$ such that both $G\backslash S$ and $G\backslash N_G[S]$ are $k$-shellable, then $G$ is $k$-shellable.
\end{thm}
\begin{proof}
Let $G'=G\backslash N_G[S]$ and $G''=G\backslash S$. Let $F_1,\ldots,F_r$ and $H_{r+1},\ldots,H_s$ be, respectively, the $k$-shelling orders of $\Del_{G''}$ and $\Del_{G'}$. Set $F_i=H_i\cup S$, for $r+1\leq i\leq s$. We show that $F_1,\ldots,F_s$ is a $k$-shelling order of $\Del_G$.

Note that $\Del_{G''}=\dl_{\Del_G}(S)$. This follows from the fact that $S$ is a $k$-simplicial set. Hence $F_1,\ldots,F_s$ contains all of facets of $\Del_G$. Now let $1\leq j<i\leq s$. If $i\leq r$ or $j>r$, by $k$-shellability of $\Del_{G'}$ and $\Del_{G''}$, we are done. So suppose that $i>r$ and $j\leq r$. Clearly, $S\subseteq F_i\backslash F_j$. On the other hand, $\Del_{G'}\subset\Del_{G''}$ and so there exists $l\leq r$ such that $H_i\subset F_l$. This implies that $S=F_i\backslash F_l$ and therefore the assertion is completed.
\end{proof}

The following theorem extends Theorem 2.1.13 of \cite{CrEs}.

\begin{thm}\label{sh iff lk}
Let $S_1$ be a $k$-simplicial set of $G$ and let $S_1,\ldots,S_r$ be the parts of  $G_{N_G[S_1]}$ and $G'_i=G\backslash N_G[S_i]$ for all $i=1,\ldots,r$. Then $G$ is $k$-shellable if and only if $G'_i$ is $k$-shellable for all $i=1,\ldots,r$.
\end{thm}
\begin{proof}
Only if part follows from Theorem \ref{G to lk}. Conversely, let $G'_i$ be $k$-shellable for all $i=1,\ldots,r$. Hence for every $i$, there exists a $k$-shelling order $F_{i1},\ldots,F_{it_i}$ for $\Del_{G'_i}=\lk_{\Del_G}(S_i)$. We claim that
$$F_{11}\cup S_1,\ldots,F_{1t_1}\cup S_1,\ldots,F_{r1}\cup S_r,\ldots,F_{rt_r}\cup S_r$$
is a $k$-shelling order for $\Del_G$.

We first show that the above list is the complete list of facets of $\Del_G$. Let $F\in\Del_G$. If for some $i\neq 1$, $S_i\subset F$ then $F$ is in the above list. Otherwise, suppose that for all $i\neq 1$, $S_i\cap F=\emptyset$. Since the elements of $S_1$ are only adjacent to the elements of $S_i$'s, for $i\neq 1$, it follows that $S_1\subset F$. On the other hand, it is not possible that $F$ contains only some of elements of a part, as $S_i$, but not all of elements of $S_i$. Because in this case, there are $x_j\in S_i\cap F$ and $x_{j'}\in S_i\backslash F$. This means that $N_G(x_j)\cap F=\emptyset$ but $N_G(x_{j'})\cap F\neq\emptyset$ and so $N_G(x_j)\neq N_G(x_{j'})$, a contradiction.

Now consider $F_{ij}\cup S_i$ and $F_{lm}\cup S_l$. We have the following cases:

(i) $i<l$: Then $S_l\subset (F_{lm}\cup S_l)\backslash (F_{ij}\cup S_i)$. Since $F_{lm}\cup S_1$ is an independence set in $G$, there exists a facet $F\in\Del_G$ with $F_{lm}\cup S_1\subset F$. In particular, $F$ is in the form $F=F_{1p}\cup S_1$ for some $1\leq p\leq t_1$. Thus $(F_{lm}\cup S_l)\backslash (F_{1p}\cup S_1)=S_l$.

(ii) $i=l$ and $j<m$: The assertion follows from $k$-shellability of $\Del_{G'_i}$.

Therefore $G$ is $k$-shellable.
\end{proof}

\begin{lem}\label{neib facet}
Let $S,T\subset V(G)$ and let $T$ be an independent set of $G$. If $N_G[T]\subseteq N_G[S]$ then independent sets of $G\backslash N_G[S]$ are not maximal independent sets of $G\backslash S$.
\end{lem}
\begin{proof}
If $F$ is an independent set of $G\backslash N_G[S]$ then $F\cup T$ will be a larger independent set of $G\backslash S$.
\end{proof}

\begin{cor}
Let $S$ be a $k$-simplicial set of $G$ and let $S_1,\ldots,S_r$ be the parts of  $G_{N_G[S]}$. Then for every $S_i$, independent sets of $G\backslash N_G[S_i]$ are not maximal independent sets of $G\backslash S_i$.
\end{cor}
\begin{proof}
Since $N_G[S]\subseteq N_G[S_i]$ for all $i$, the assertion follows from Lemma \ref{neib facet}.
\end{proof}

In \cite[Theorem 2]{CaCr} the authors proved that if $x$ is a simplicial vertex of $G$ and $y$ adjacent to $x$, then $G$ is shellable if and only if $G\backslash N_G[y]$ and $G\backslash y$ are shellable. The following theorem extends this result to $k$-shellable graphs.

\begin{thm}\label{k-sh equi}
Let $S_1$ be a $k$-simplicial set of $G$ and let $S_1,\ldots,S_r$ be the parts of  $G_{N_G[S_1]}$. Then $G$ is $k$-shellable if and only if for each $i=2,\ldots,r$ the graphs $G'_i=G\backslash N_G[S_i]$ and $G''_i=G\backslash S_i$ are $k$-shellable.
\end{thm}
\begin{proof}
``Only if part'': Let $G$ be $k$-shellable. It follows from Theorem \ref{sh iff lk} that $G'_i$ is $k$-shellable for all $i$. Fix an integer $i$. We want to show that $G''_i$ is $k$-shellable. By again relabeling $S_i$'s we can consider $i=r$.

Since each $G'_i$ is $k$-shellable, so for every $i$, there exists a $k$-shelling order $F_{i1},\ldots,F_{it_i}$ for $\Del_{G'_i}=\lk_{\Del_G}(S_i)$. Moreover, by the proof of Theorem \ref{sh iff lk},
$$F_{11}\cup S_1,\ldots,F_{1t_1}\cup S_1,\ldots,F_{r1}\cup S_r,\ldots,F_{rt_r}\cup S_r$$
is a $k$-shelling order for $\Del_G$. By the fact that $N_G(S_1)\subseteq N_G(S_r)$ we conclude that for every $F_{rj}$ where $1\leq j\leq t_r$, there exists $1\leq l\leq t_1$ such that $F_{rj}\subseteq F_{1l}$. Therefore
$$F_{11}\cup S_1,\ldots,F_{1t_1}\cup S_1,\ldots,F_{r-1\ 1}\cup S_{r-1},\ldots,F_{r-1\ t_{r-1}}\cup S_{r-1}$$
will be a list of facets of $G\backslash S_r$. Furthermore, it is a $k$-shelling order of $G''_r$.

``If part'': Let for all $i=2,\ldots,r$ the graphs $G'_i$ and $G''_i$ are $k$-shellable. Fix an $i$ and set $G'=G'_i,G''=G''_i,S=S_i$. Let $F_1,\ldots,F_r$ and $H_1,\ldots,H_s$ be, respectively, the $k$-shelling orders of $\Del_{G''}$ and $\Del_{G'}$. We first show that $$F_1,\ldots,F_r,H_1\cup S,\ldots,H_s\cup S$$ is a list of facets of $\Del_G$, and furthermore, this list is a $k$-shelling order of $\Del_G$.

Let $F$ be a facet of $\Del_G$. If $F\cap S=\emptyset$, then since $F$ contains no edge of $G$, it contains no edge of $G''$ and so $F\in\Del_{G''}$. Suppose $F\cap S\neq\emptyset$. Let $x\in F\cap S$ and let $y\in S$ with $y\neq x$. Since $N_G(y)=N_G(x)$, we have $y\in F$. Thus $S\subset F$. On the other $F$ contains no edge of $G'$. Therefore $F\backslash S$ is a facet of $\Del_{G'}$.

Now we show that above list is a $k$-shelling order. Set $F_{i+r}=H_i\cup S$ for all $i=1,\ldots,s$. Suppose $F_i$ and $F_j$ with $i<j$. If $j\leq r$ or $i\geq r+1$ then by the $k$-shellability of $\Del_{G''}$ and $\Del_{G'}$, respectively, the assertion is completed. Let $i\leq r$ and $r<j$. Then $S\subset F_j\backslash F_i$. Since $H_{j-r}\in\Del_{G''}$, there is $F_l$ with $l\leq r$ such that $H_{j-r}\subset F_l$. Therefore $F_j\backslash F_l=S$. This completes the assertion.
\end{proof}

\textbf{Acknowledgments:}
The author would like to thank the referee for careful reading of the paper. The work was supported by the research council of the University of Maragheh.

\end{document}